\documentclass[a4paper,12pt]{article}
\usepackage{graphicx}
\usepackage{textcomp,graphicx,epsfig}
\usepackage{latexsym,amssymb,amsmath}

\newcommand{\PP}{\mathbb{P}}

\newcommand{\ee}{\mathrm{e}}
\newcommand{\ii}{\mathrm{i}}

\newtheorem{defi}{Definition}[section]

\newtheorem{prop}[defi]{Proposition}

\newenvironment{dem}{\vskip 2mm\noindent {\it Proof}:}
                    {\hfill $\square$ \vskip 2mm \noindent}

\begin{document}

\title{Computing wedge probabilities
\thanks{Research supported by laboratoire d'excellence TOUCAN (Toulouse Cancer)}
}

\author{Bernard Ycart and R\'emy Drouilhet\footnote{
          Laboratoire Jean Kuntzmann,
          Universit\'e Grenoble Alpes,
          51 rue des math\'ematiques, 38041 Grenoble cedex, France
          \texttt{bernard.ycart@imag.fr}
}}

\maketitle

\begin{abstract}
A new formula for the probability that a standard Brownian motion stays between 
two linear boundaries is proved. A simple algorithm is deduced. Uniform
precision estimates are computed. 
Different implementations have been made available online as
\texttt{R} packages.
\end{abstract}

\section{Introduction}
\label{intro}
Let $W=\{W_t\,,\;t\geqslant 0\}$ be a standard Brownian motion
defined on a filtered probability space $(\Omega,\mathcal{F},\PP)$.
The probability that $W_t$ remains in the planar region between
two linear boundaries $-a_1 t -b_1$ and $a_2 t +b_2$ 
will be referred to as \emph{wedge probability}, and denoted by
$k(a_1,b_1;a_2,b_2)$:
\begin{equation}
\label{wedge}
k(a_1,b_1;a_2,b_2)
=
\PP[-a_1 t -b_1 \leqslant W_t \leqslant a_2 t + b_2
\,,\;\mbox{for all } t\geqslant 0]\;.
\end{equation}
It is positive if and only if 
$a_1,b_1,a_2,b_2$ are all positive, which will be assumed
from now on. Doob \cite{Doob49} expressed $k(a_1,b_1;a_2,b_2)$ as the sum of
a convergent series. Since then, Doob's formula has been extended or
applied by many authors, including 
\cite{Anderson60,Durbin71,Kennedy73,BarbaEscriba87,Hall97,SalminenYor11}.
One reason for its success is that many boundary crossing
problems reduce to computing a wedge probability, through the
representation of a certain Gaussian process in terms of $W$
\cite{Park81,BorodinSalminen02}. The earliest
example is the standard Brownian bridge:
$$
\left\{B_{t},0\leqslant t<1\right\}
\overset{d}{=}
\left\{(1-t)W_{t/(1-t)},0\leqslant t<1\right\}\,.
$$%
From this representation one gets:
\begin{equation}
\label{KS}
\PP[\sup_{t\in[0,1]} |B_t|\leqslant a] = k(a,a;a,a)\;,
\end{equation}
which is the 
the distribution function of the test statistic in the 
Kolmogorov-Smirnov two-sided test, computed by Kolmogorov \cite{Kolmogorov33};
see \cite{Stephens92} for historical aspects. More generally, 
the probability that a (non necessarily standard) Brownian bridge
stays between two linear segments is a wedge probability. This remark 
makes wedge probabilities a building block for more general
boundary crossing problems, through the method of piecewise linear 
approximations.
\vskip 1mm\noindent
The exit probability of a stochastic process from a region of the
plane limited by two curves is called 
Boundary Crossing Probability (BCP). Applications of BCP's
can be found in many fields, from non-parametric
statistics to biology or finance: 
see \cite{WangPotzelberger07} and references therein.
Explicit results are scarse
\cite{Kahale08}. A general approximation method has been proposed by
Wang and Potzelberger \cite{WangPotzelberger97} for single boundaries and
Novikov et al. \cite{Novikovetal99} for double boundaries; see also
\cite{PotzelbergerWang01,BorovkovNovikov05,WangPotzelberger07}.
The idea is to replace the two (nonlinear) boundaries by piecewise linear
approximations. Given its two values at the bounds of an interval, the
conditional distribution of $W$ is that of a Brownian bridge on that
interval. Thus the probability that it stays between two linear
segments is a wedge probability. Using the independent increment
property of $W$,
the probability that the standard Brownian motion
stays between two piecewise linear boundaries is written as a
multidimensional Gaussian integral, the integrand being a product of
wedge probabilities \cite[Theorem 1]{Novikovetal99}.
The integral can be approximated either as a Gauss-Hermite quadrature 
\cite{GenzBretz09} or by a Monte Carlo method \cite{Potzelberger12}.
In both cases, the integrand must be repeatedly evaluated, 
which implies that many wedge
probabilities must be calculated for very different sets of
values. The problem is that in Doob's formula, as well as in all other
equivalent formulas published since, the speed of convergence of the
series depends on the parameters, and may be very slow for small
values. This makes the BCP approximation algorithms numerically
unstable.
\vskip 1mm\noindent
The key to efficient computation of wedge probabilities has long been 
available: Jacobi's theta functions and their double expression
through Poisson's summation formula. Kolmogorov \cite{Kolmogorov33} had already 
given two formulas for $k(a,a;a,a)$, and remarked the interest for
numerical computation: one converges fast for relatively large values
of $a$, the other for relatively small values. This is routinely used
in statistical softwares implementing the Kolmogorov-Smirnov test. The
connection of $k(a_1,b_1;a_2,b_2)$ with theta functions has been
pointed out by Salminen and Yor 
\cite{SalminenYor11}. However, no alternative to Doob's
formula has been deduced so far; this is the main contribution of this
paper (Proposition \ref{prop:wedge}). The algorithmic consequence is
that computing at most three terms of the series either in 
Doob's formula or in the new
alternative suffices to approximate $k(a_1,b_1;a_2,b_2)$ with
precision smaller than $10^{-16}$. Uniform bounds on precision and
algorithmic consequences are presented in section 
\ref{sec:algo}. Several implementations are compared 
in section \ref{sec:implementation}.
An \texttt{R} package \texttt{wedge} 
has been made available online
\cite{DrouilhetYcart16}. Its companion \texttt{wedgeParallel}
permits full use of a multicore structure.
\section{Alternative to Doob's formula}
\label{sec:wedge}
Doob \cite{Doob49} (formula (4.3) p.~398) expresses 
$k(a_1,b_1;a_2,b_2)$ as follows:
\begin{equation}
\label{doob}
k(a_1,b_1;a_2,b_2)
= 1 - \sum_{n=1}^{+\infty} 
\ee^{-2A_n}+\ee^{-2B_n}
-\ee^{-2C_n}-\ee^{-2D_n}\;,
\end{equation}
with:
\begin{eqnarray*}
A_n &=& n^2a_2b_2+(n-1)^2a_1b_1
+n(n-1)(a_2b_1+a_1b_2)\;,\\
B_n &=& (n-1)^2a_2b_2+n^2a_1b_1
+n(n-1)(a_2b_1+a_1b_2)\;,\\
C_n &=& n^2(a_1b_1+a_2b_2)
+n(n-1)a_2b_1+n(n+1)a_1b_2\;,\\
D_n &=& n^2(a_1b_1+a_2b_2)
+n(n+1)a_2b_1+n(n-1)a_1b_2\;.
\end{eqnarray*}
Here is another expression.
\begin{prop}
\label{prop:wedge}
For $a_1,b_1,a_2,b_2>0$, denote:
$$
a_\pm = \scriptstyle{\frac{a_1\pm a_2}{2}}
\displaystyle{\,;\;
b_\pm =}\scriptstyle{\frac{b_1\pm b_2}{2}}
\displaystyle{\,;\;
c = }\scriptstyle{\frac{a_1b_1 - a_2b_2}{2}}
\displaystyle{\,;\;
d = }\scriptstyle{\frac{a_1b_2 - a_2b_1}{2}}
\displaystyle{\;.}
$$
Then:
\begin{equation}
\label{wedge2}
\begin{array}{rl}
\displaystyle{
k(a_1,b_1;a_2,b_2)=}
\scriptstyle{\sqrt{\frac{\pi}{2a_+b_+}}}
\displaystyle{\,
\ee^{\frac{d^2}{2a_+b_+}}
\sum_{n=1}^{+\infty}}&\displaystyle{
\Big(\ee^{-\frac{\pi^2(2n)^2}{8a_+b_+}}
\big(\cos( \scriptstyle{\frac{\pi (2n)d}{2a_+b_+}}
\displaystyle{)-\cos(}
\scriptstyle{\frac{\pi (2n)c}{2a_+b_+}}\displaystyle{)\big)}}
\\
&\displaystyle{+
\ee^{-\frac{\pi^2(2n-1)^2}{8a_+b_+}}
\big(\cos(\scriptstyle{\frac{\pi (2n-1)d}{2a_+b_+}}
\displaystyle{)+\cos(}
\scriptstyle{\frac{\pi (2n-1)c}{2a_+b_+}}
\displaystyle{)\big)\Big)\;.}}
\end{array}
\end{equation}
\end{prop}
\begin{dem}
In terms of $a_\pm$, $b_\pm$, $c$, $d$,
the expressions of $A_n$, $B_n$, $C_n$, $D_n$ are:
\begin{eqnarray*}
A_n &=&
(2n-1)^2 a_+b_+
-(2n-1) c
+ a_-b_-
\;,\\
B_n &=& 
(2n-1)^2 a_+b_+
+(2n-1) c
+ a_-b_-
\;,\\
C_n &=& (2n)^2 a_+b_+ + 2n d\;,
\\
D_n &=& (2n)^2 a_+b_+ - 2n d\;.
\end{eqnarray*}
Hence:
\begin{equation}
\label{wedge1}
\begin{array}{rl}
\displaystyle{
k(a_1,b_1;a_2,b_2)=1+2\sum_{n=1}^{+\infty}}
&\displaystyle{\Big( 
\ee^{-2(2n)^2a_+b_+}\cosh(2(2n)d)}
\\
&\displaystyle{
-\ee^{-2a_-b_-} \ee^{-2(2n-1)^2a_+b_+}
\cosh(2(2n-1)c)\Big)\;.}
\end{array}
\end{equation}
Not meaning to add anything to the ``bewildering variety of notations''
for theta functions \cite[p.~576]{AbramowitzStegun64}, let us denote by
$\theta$ the following function of two complex variables:
\begin{equation}
\label{deftheta}
\theta(u,v) = \sum_{n=-\infty}^{+\infty} \ee^{-2n^2 u} \cosh(2n v)\;. 
\end{equation}
Observe that:
\begin{equation}
\label{defthetai}
\theta(u,v+\scriptstyle{\frac{\ii \pi}{2}}
\displaystyle{)} = \sum_{n=-\infty}^{+\infty} (-1)^n \ee^{-2n^2 u} \cosh(2n v)\;. 
\end{equation}
By Poisson's summation formula  (see for instance 
formula (11) p.~236 of
\cite{VanderPolBremmer50}), one gets:
\begin{equation}
\label{defthetainv}
\theta(u,v) = \sqrt{\frac{\pi}{2u}}\,\ee^{v^2/(2u)}
\sum_{n=-\infty}^{+\infty} \ee^{-\pi^2n^2/(2u)}\cos(\pi n v/u)\;, 
\end{equation}
and:
\begin{equation}
\label{defthetaiinv}
\theta(u,v+\scriptstyle{\frac{\ii \pi}{2}}
\displaystyle{)}
= \sqrt{\frac{\pi}{2u}}\,\ee^{v^2/(2u)}
\sum_{n=-\infty}^{+\infty} \ee^{-\pi^2(n+\frac{1}{2})^2/(2u)}
\cos(\pi (n+\scriptstyle{\frac{1}{2}}\displaystyle{) v/u)\;.} 
\end{equation}
From (\ref{wedge1}), (\ref{deftheta}), and (\ref{defthetai}):
\begin{eqnarray*}
k(a_1,b_1;a_2,b_2)
&=&
\frac{1}{2}
\big(
\theta(a_+b_+,d)+
\theta(a_+b_+,d+\scriptstyle{\frac{\ii \pi}{2}}
\displaystyle{)}
\big)\\
&&
-
\frac{\ee^{-2a_-b_-}}{2}
\big(
\theta(a_+b_+,c)-
\theta(a_+b_+,c+\scriptstyle{\frac{\ii \pi}{2}}
\displaystyle{)}
\big)\;.
\end{eqnarray*}
\vskip 1mm\noindent
Combining four evaluations of $\theta$
does not quite solve the numerical problem for small values of
$a_+b_+$. The
terms of the four series need further rearrangement. It is obtained
observing that: 
$$
c^2-4a_-b_-a_+b_+ = d^2\;.
$$
Hence:
$$
\ee^{-2a_-b_-}\ee^{\frac{c^2}{2a_+b_+}} = 
\ee^{\frac{d^2}{2a_+b_+}}\;.
$$
From there, (\ref{wedge2}) follows.
\end{dem}
\vskip 1mm\noindent
As remarked by Salminen and Yor 
\cite{SalminenYor11}, the symmetry and scaling properties 
of $k$ can be read on Doob's formula. They also appear on
(\ref{wedge2}):
$$
k(a_1,b_1;a_2,b_2) = k(a_2,b_2;a_1,b_1)=k(b_1,a_1;b_2,a_2)
=k(\frac{a_1}{u},ub_1;\frac{a_2}{u},ub_2)\;.
$$
When slopes or intercepts
are equal the expressions are simpler. If $a_1=a_2=a$ then $a_+=a$, $a_-=0$, 
$c=-d=a b_-$, and:
\begin{eqnarray*}
k(a,b_1;a,b_2)&=&
1+2\sum_{n=1}^{+\infty} (-1)^n
\ee^{-2n^2a  b_+}\cosh(2nab_-)\\
&=&
\sqrt{\frac{2\pi}{ab_+}}
\ee^{\frac{ab_-^2}{2b_+}}
\sum_{n=1}^{+\infty} \ee^{-\frac{\pi^2(2n-1)^2}{8ab_+}}
\cos(\pi(2n-1)\scriptstyle{\frac{b_-}{2b_+}}
\displaystyle{)) \;.}
\end{eqnarray*}
When both slopes and intercepts are equal one gets:
\begin{eqnarray*}
k(a,b;a,b)&=&
1+2\sum_{n=1}^{+\infty} (-1)^{n}
\ee^{-2n^2ab}\\
&=&
\sqrt{\frac{2\pi}{ab}}
\sum_{n=1}^{+\infty} \ee^{-\frac{\pi^2(2n-1)^2}{8ab}}
\;.
\end{eqnarray*}
The first sum is formula (4.3') of \cite{Doob49}.
\vskip 1mm\noindent
When slopes equal intercepts, the probability for a standard Brownian
bridge to stay  
in a horizontal band is obtained, i.e. formula (4.9) p.~448 of \cite{Bianeetal01}.
If $b_1=a_1$ and $b_2=a_2$, then $b_\pm=a_\pm$,
$c=2a_+a_-$, $d=0$, and:
\begin{eqnarray*}
k(a_1,a_1;a_2,a_2)&=&
1+\sum_{n=1}^{+\infty} 2\ee^{-2(n(a_1+a_2))^2}
-\ee^{-(2n(a_1+a_2)-a_1)^2}
-\ee^{-(2n(a_1+a_2)-a_2)^2}
\\
&=&
\sqrt{\frac{\pi}{2}}\frac{1}{a_+}
\sum_{n=1}^{+\infty} 
\ee^{-\frac{\pi^2(2n)^2}{8a_+^2}}
(1-\cos(\pi(2n)
\scriptstyle{\frac{a_-}{a_+}}\displaystyle{))
}\\
&&\hspace*{1.6cm}+
\ee^{-\frac{\pi^2(2n-1)^2}{8a_+^2}}
(1+\cos(\pi(2n-1)
\scriptstyle{\frac{a_-}{a_+}}\displaystyle{))
\;.}
\end{eqnarray*}
\vskip 1mm\noindent
Finally, the case
where all four parameters are equal is the 
probability for a standard Brownian bridge to stay
in a horizontal band centered at $0$, i.e. 
the distribution function of the test statistic in the 
Kolmogorov-Smirnov two-sided test.
The formulas were originally found by Kolmogorov
\cite{Kolmogorov33}; Feller
\cite{Feller48} gave a simpler proof. 
Both had remarked the double expression coming from theta
functions, and its interest for numerical computation.
\begin{eqnarray*}
k(a,a;a,a)&=&
1+2\sum_{n=1}^{+\infty} (-1)^{n}
\ee^{-2n^2 a^2}\\
&=&
\frac{\sqrt{2\pi}}{a}\sum_{n=1}^{+\infty} \ee^{-\frac{\pi^2(2n-1)^2}{8a^2}}\;.
\end{eqnarray*}
\section{Algorithm and precision}
\label{sec:algo}
\begin{figure}[!ht]
\centerline{
\includegraphics[width=12cm]{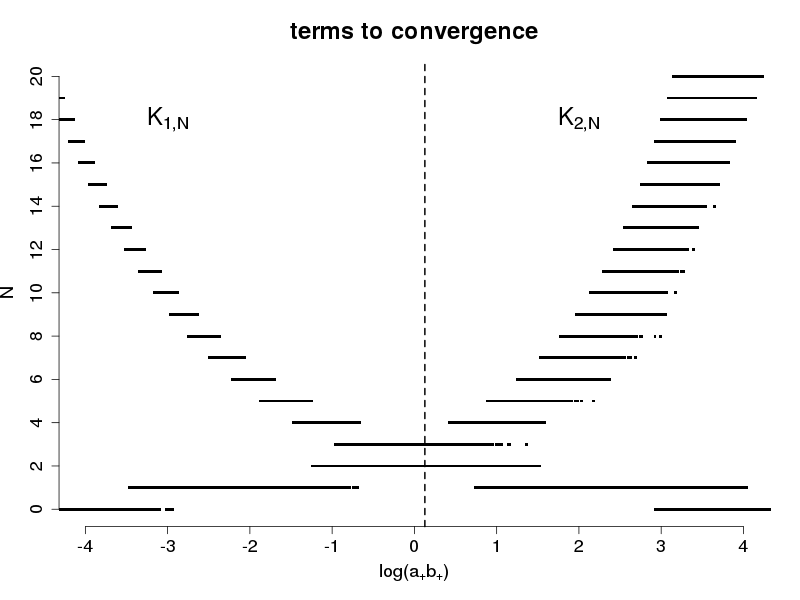}
} 
\caption{Number of terms to convergence as a function of $\log(a_+b_+)$ 
for $K_{1,N}$ and $K_{2,N}$, over $10^6$ simulated values of 
$a_1,b_1;a_2,b_2$.  The dashed vertical line marks the theoretical
threshold for $N=3$ i.e. $\tau_3=1.136$.}
\label{fig:convabpbw}
\end{figure}

Denote by $K_{1,N}$ and
$K_{2,N}$ the partial sums up to $N$ in formulas (\ref{doob}) and 
(\ref{wedge2}).
\begin{equation}
\label{doobN}
K_{1,N}
= 1 - \sum_{n=1}^{N} 
\ee^{-2A_n}+\ee^{-2B_n}
-\ee^{-2C_n}-\ee^{-2D_n}\;,
\end{equation}
\begin{equation}
\label{wedge2N}
\begin{array}{rl}
\displaystyle{K_{2,N} = 
\sqrt{\frac{\pi}{2a_+b_+}}\,
\ee^{\frac{d^2}{2a_+b_+}}
\sum_{n=1}^{N}}&\displaystyle{
\Big(\ee^{-\frac{\pi^2(2n)^2}{8a_+b_+}}
\big(\cos( \scriptstyle{\frac{\pi (2n)d}{2a_+b_+}}
\displaystyle{)-\cos(}
\scriptstyle{\frac{\pi (2n)c}{2a_+b_+}}\displaystyle{)\big)}}\\
&\displaystyle{+
\ee^{-\frac{\pi^2(2n-1)^2}{8a_+b_+}}
\big(\cos(\scriptstyle{\frac{\pi (2n-1)d}{2a_+b_+}}
\displaystyle{)+\cos(}
\scriptstyle{\frac{\pi (2n-1)c}{2a_+b_+}}
\displaystyle{)\big)\Big)\;.}}
\end{array}
\end{equation}
The question is: for a given set of parameters
$a_1,b_1,a_2,b_2$, which of $K_{1,N}$ 
and $K_{2,N}$ should be computed,
and which value of $N$ ensures a given precision?
Proposition \ref{prop:remainders}
bounds remainders. 
\begin{prop}
\label{prop:remainders}
Denote by $R_{1,N}$ and $R_{2,N}$ the remainders:
$$
R_{1,N} = K_{1,\infty}-K_{1,N}
\quad\mbox{and}\quad
R_{2,N} = K_{2,\infty}-K_{2,N}\;.
$$
For $N>1$:
\begin{equation}
\label{rem1N}
R_{1,N} \leqslant 
\scriptstyle{\frac{1}{4a_+b_+(N-1)}}
\displaystyle{
\ee^{-8a_+b_+(N-1)^2}\;,}
\end{equation}
\begin{equation}
\label{rem2N}
R_{2,N} \leqslant
\scriptstyle{
\left(\frac{2}{\pi}\right)^{3/2}
\frac{\sqrt{a_+b_+}}{N}}
\displaystyle{ \ee^{2a_+b_+} 
\ee^{-\frac{\pi^2 N^2}{2a_+b_+}}
\;.}
\end{equation}
\end{prop}
\begin{dem}
For both bounds, the following well known 
inequality is used: for any positive $u$,
\begin{equation}
\label{majsumn2}
\sum_{n=N+1}^{+\infty} \ee^{-un^2}
\leqslant
\frac{\ee^{-uN^2}}{2uN}\;.
\end{equation}
%
In $R_{1,N}$, all four terms $A_n$, $B_n$, $C_n$, $D_n$ are larger than
$4(n-1)^2a_+b_+$. Hence:
$$
R_{1,N} \leqslant
4 \sum_{n=N+1}^{+\infty}\ee^{-8(n-1)^2a_+b_+}\;,
$$
from where (\ref{rem1N}) follows by (\ref{majsumn2}).
\vskip 1mm\noindent
For $R_{2,N}$, notice first that
$|d|\leqslant
2\sqrt{a_+b_+}$, by Schwarz inequality.
Hence:
\begin{eqnarray*}
R_{2,N} &\leqslant&
\sqrt{\frac{2\pi}{a_+b_+}} \ee^{2a_+b_+} \sum_{n=N+1}^{+\infty}
\ee^{-\frac{\pi^2(2n)^2}{8a_+b_+}}
+
\ee^{-\frac{\pi^2(2n-1)^2}{8a_+b_+}}\\
&=&
\sqrt{\frac{2\pi}{a_+b_+}} \ee^{2a_+b_+}
\sum_{n=2N+1}^{+\infty}\ee^{-\frac{\pi^2n^2}{8a_+b_+}}\;.
\end{eqnarray*}
Using again (\ref{majsumn2}) leads to (\ref{rem2N}).
\end{dem}
As expected, the bound on $R_{1,N}$ decreases with $a_+b_+$, 
the bound on $R_{2,N}$ increases. 
Denote by $\tau_N$ the value of $a_+b_+$ such that both bounds are equal, and by 
$\varepsilon_N$ their common value. Computing 
$K_{1,N}$ if $a_+b_+\geqslant \tau_N$ and $K_{2,N}$ else 
ensures $\varepsilon_N$ precision at least, whatever $a_1,b_1,a_2,b_2$. 
The values of $\tau_N$ and 
$\varepsilon_N$ for $N=2,\ldots,8$ are given in Table \ref{tab:epsN}.
\begin{table}[!ht]
\centerline{
\begin{tabular}{c||c|c|c|c|c|c|c}
$N$&$2$&$3$&$4$&$5$&$6$&$7$&$8$\\
\hline
$\tau_N$&$1.380$&$1.136$&$1.030$&$0.973$&$0.937$&$0.912$&$0.895$\\
$\varepsilon_N$&$\scriptstyle{2.9~10^{-6}}$&$\scriptstyle{1.8~10^{-17}}$
&$\scriptstyle{5.1~10^{-34}}$&$\scriptstyle{5.6~10^{-56}}$
&$\scriptstyle{2.3~10^{-83}}$&$\scriptstyle{3.5~10^{-116}}$
&$\scriptstyle{1.9~10^{-154}}$
\\
\hline
\end{tabular}
}
\caption{Threshold and precision per number of terms computed.}
\label{tab:epsN}
\end{table}
In particular, for $N=3$ a precision $\varepsilon_3=1.8~10^{-17}$ is obtained,
which is below current machine double precision. Thus $N=3$ was chosen
as default value 
in our implementation.
\vskip 1mm\noindent
Computing more than two terms is usually not necessary. To illustrate this
a Monte Carlo study has been conducted, over $10^6$ four-tuples of 
independent random values
drawn in the interval $[0{,}10]$ with cumulative distribution function 
$(x/10)^{1/2}$. That choice ensured that wedge probabilities covered 
the whole range of $[0,1]$, with higher mass on 
values close to $0$ or $1$. For both sums the number of terms to convergence 
was defined as the first value of $N$ such that the remainder is smaller
than $\varepsilon=10^{-16}$. As predicted by Proposition \ref{prop:remainders},
in all $10^6$ cases either $K_{1,N}$ or $K_{2,N}$ reached $\varepsilon$ precision 
with $N=3$ terms or less. Actually, 
in $2.5\%$ of the cases, the result was smaller than $\varepsilon$ or larger
than $1-\varepsilon$: no summation was needed. In $73.8\%$ of the cases 
$N=1$ sufficed to get $\varepsilon$ precision, and in $20.6\%$ of the cases 
$N=2$ terms were necessary; only in $0.56\%$ of the cases were $N=3$ terms
necessary. 
Experimental results evidenced the need for 
an alternative to $K_{1,N}$. Indeed, in $442$ out of the $10^6$
cases, the number of terms to convergence of $K_{1,N}$ 
was larger than $100$, and in $1558$ cases
it was larger than $50$.
Figure \ref{fig:convabpbw} presents the numbers of terms to convergence
as a function of $\log(a_+b_+)$, for all $10^6$ random values, and both
$K_{1,N}$ and $K_{2,N}$. As expected, the numbers decrease for 
$K_{1,N}$; they increase for $K_{2,N}$.
\section{Implementation}
\label{sec:implementation}
The calculation of successive terms in
$K_{1,N}$ or $K_{2,N}$ is easily vectorized. This makes the
computation of a vector of wedge probabilities relatively fast in pure
\texttt{R} \cite{R}.
Our objective was to explore gains in computing time,
using existing \texttt{R} tools. The most widely 
used of these tools is \texttt{Rcpp} \cite{Rcpp}. 
It uses (usually faster) compiled 
\verb@C++@ code, interfaced with the \texttt{R} environment.
Most computers now have a multicore architecture. 
However by default, both \texttt{R} or \texttt{Rcpp} use only one core. 
Taking full advantage of a multicore stucture 
can be done for example through 
\texttt{RcppParallel} \cite{RcppParallel}.

Numerical experiments have been made using vectors of simulated entries
with the same distribution as in section \ref{sec:algo}: independent entries
on $[0{,}10]$ with cumulative distribution function 
$(x/10)^{1/2}$. Five implementations were considered: pure
\texttt{R}, 
one-core 
\texttt{Rcpp}, \texttt{RcppParallel} with 4, 6, and 8 cores. 
Table 2 reports
running times on a MacBookPro Retina 15. The running time for $10^6$
values in pure \texttt{R} 
($0.725$ second), can be considered satisfactory. 
However, the gain in time
goes up to twentyfold if an eight-core architecture is used. One of the
known limitations of 
vectorized versions in pure \texttt{R} 
is memory space: ours cannot deal with vectors
larger than $10^7$ entries.

\begin{table}[!ht]
\label{tab:resexpe}
\centerline{
\begin{tabular}{c||c|c|c|c|c}
$n$ & \texttt{R} & \texttt{Rcpp} & \texttt{RcppParallel (4)} & \texttt{RcppParallel (6)} & \texttt{RcppParallel (8)} \\\hline
$10^6$ & 0.725 & 0.169 & 0.05 & 0.044 & 0.036 \\
$10^7$ & 5.854 & 1.785 & 0.471 & 0.417 & 0.353 \\
$10^8$ & --- & 17.607 & 4.725 & 4.266 & 3.987 \\\hline
\end{tabular}
}
\caption{Times in second for calculating wedge probabilities over 
vectors of size $n=10^6,10^7,10^8$, in pure \texttt{R},
\texttt{Rcpp}, \texttt{RcppParallel} with 4, 6, 8 
threads. For $n=10^8$, calculations exceed memory space in pure \texttt{R}.}
\end{table}

An \texttt{R} 
package \texttt{wedge} has been made available online \cite{DrouilhetYcart16}. 
In order to address installation issues for users not interested by a
parallel version, the companion package
\texttt{wedgeParallel} has been left as an option.


\bibliographystyle{plain}

%
%

\end{document}